\documentclass[a4paper,12pt]{amsart}
\usepackage{amssymb}
\usepackage{amsmath}

\def\a{\alpha}

\def\vs{\vskip .6cm}

\def\beq{\begin{equation}}
\def\eeq{\end{equation}}
\def\bea{\begin{eqnarray*}}
\def\eea{\end{eqnarray*}}
\def\beaa{\begin{eqnarray}}
\def\eeaa{\end{eqnarray}}
\def\ba{\begin{array}}
\def\ea{\end{array}}

\def\e{\varepsilon}

\def \RM{\mathbb{R}}

\def \ZM{\mathbb{Z}}
\def \CM{\mathbb{C}}

\def \HM{\mathbb{H}}

\def\rank{\mathrm{rk}}
\def\i{\mathrm{i}}
\def\rk{\mathrm{rk}}

\def\id{\mathrm{id}}
\def\be{\begin{equation}}
\def\ee{\end{equation}}

\def\rk{\mathrm{rk}}

\def\Sym{\mathrm{Sym}}

\def\gg{\mathfrak{g}}

\def\tt{\mathfrak{t}}

\def\kk{\kappa}
\def\mm{\mathfrak{m}}

\def\SU{\mathrm{SU}}
\def\U{\mathrm{U}}
\def\E{\mathrm{E}}
\def\H{\mathrm{H}}
\def\R{\mathrm{R}}
\def\F{\mathrm{F}}
\def\SO{\mathrm{SO}}

\def\End{\mathrm{End}}

\def\Sp{\mathrm{Sp}}
\def\Spin{\mathrm{Spin}}

\def\Sym{\mathrm{Sym}}

\def\ind{\mathrm{ind}}
\def\Gr{\mathrm{Gr}}
\def\ch{\mathrm{ch}}

\def\P{\mathrm{P}}
\def\T{\mathrm{T}}
\def\S{\mathrm{S}}
\def\D{\mathrm{D}}
\def\bb{\mathrm{b}}
\def\RR{{\mathcal R}}

\newtheorem{epr}{Proposition}[section]

\newtheorem{ath}[epr]{Theorem}

\newtheorem{elem}[epr]{Lemma}

\theoremstyle{definition}

\newtheorem{ere}[epr]{Remark}

\title[Almost complex structures]{Almost complex structures on quaternion-K\"ahler manifolds and
  inner symmetric spaces}

\author{Paul Gauduchon, Andrei Moroianu, Uwe Semmelmann}

\address{Paul Gauduchon and Andrei Moroianu \\ CMLS\\ {\'E}cole
  Polytechnique \\ UMR 7640 du CNRS
\\ 91128 Palaiseau \\ France}
\email{pg@math.polytechnique.fr, am@math.polytechnique.fr}

\address{Uwe Semmelmann\\ Mathematisches Institut, Universit{\"a}t zu
K{\"o}ln\\
Weyertal 86-90 D-50931 K{\"o}ln, Germany}
\email{uwe.semmelmann@math.uni-koeln.de}

\date{\today}

\begin{document}

\begin{abstract}
We prove that compact quaternionic-K\"ahler manifolds of positive
scalar curvature admit no almost complex structure, even in the weak
sense, except for the complex Grassmannians
$\Gr_2(\CM^{n+2})$. We also prove that irreducible inner
symmetric spaces $M^{4n}$ of compact type 
are not weakly complex, except
for spheres and Hermitian symmetric spaces.

\vs


\end{abstract}

\maketitle

\section{Introduction}

It is a well-known fact  that the quaternionic projective spaces $\HM
\P^n$ have no almost complex structure.
The proof goes back to F. Hirzebruch in 1953  for $n\ge 4$
(cf. \cite{h53}). The non-existence of almost-complex structure on
$\HM \P^1=\S^4$ had been established a few
years earlier by
Ch. Ehresmann \cite{e49} and H. Hopf \cite{h49}. According to
Hirzebruch's lecture at the 1958 ICM
\cite{h58}, J. Milnor had in the meantime settled the remaining cases
$n= 2$ and $3$, but his proof has remained unpublished.
Later on, W.S.~Massey \cite{m62} gave an original proof of the
non-existence of almost-complex structure on $\HM
\P^n$, for any $n$, based on
the explicit calculation  of the ring $\mathrm{K}(X)$ and of the Chern
character $\ch(\T X)$ for $X= \HM \P^n$.

Quaternionic projective spaces are particular examples of
{\em quaternion-K\"ahler manifolds}. These, we recall, are (oriented)
$4n$-dimensional
Riemannian manifolds, whose holonomy is
contained in $\Sp(n)\cdot\Sp(1) \subset \SO (4n)$, if $n > 1$, or, if
$n = 1$,
(oriented) Einstein, self-dual $4$-dimensional Riemannian
manifolds.
In all dimensions $4n$, $n \geq 1$, quaternion-K\"ahler manifolds are
Einstein and are called {\it of positive type} if their scalar
curvature is positive. In this paper, we only consider
quaternion-K\"ahler of
positive type and we implicitly assume that they are
complete, hence compact.

For $n \geq 2$, the above definition of quaternion-K\"ahler manifolds
is equivalent to the
existence of locally defined almost complex structures $I,J, K$,
satisfying the quaternion relations and spanning a global rank 3 sub-bundle
$\mathrm{Q} \subset \End(\T M)$, which is preserved by the Levi-Civita
connection.
Almost complex structures on $M$ which are sections of
$\mathrm{Q}$  are called {\it compatible}. In
\cite{amp98}, it is shown that
quaternion-K\"ahler manifolds of positive
type admit no {\it compatible}  almost complex structure.
In particular the natural complex structure of the complex
Grassmannians $\Gr_2(\CM^{n+2})$, which constitute a well-known class
of quaternion-K\"ahler manifolds of positive type (cf. below), is not
compatible.

The first main result of this paper is:

\begin{ath}\label{main}
Let $M^{4n}$, $n\ge 2$, be a compact quaternion-K\"ahler manifold of
positive type, which is not isometric to the complex Grassmannian
$\Gr_2(\CM^{n+2})$.
Then $M ^{4n}$ has no weak almost complex structure, in the sense that the
tangent bundle $\T M$ is not stably isomorphic to a complex vector bundle.
\end{ath}

Notice that the assumption $n\ge 2$ is necessary, since $\HM \P^1=\S^4$
is weakly complex but not almost complex.

At the moment, the only known quaternion-K\"ahler manifolds of
positive type are the so-called (symmetric) {\it Wolf spaces}, namely \cite{wolf65}:
\begin{enumerate}

\item[(i)] the quaternionic projective spaces $
  \HM \P^n = \frac{\Sp(n+1)}{\Sp(n)\times \Sp(1)}$,

\item[(ii)]  the Grassmannians
  ${\Gr}_2(\CM^{n+2}) = \frac{\U(n+2)}{\U(n)\times \U(2)}$ of complex
  $2$-planes in $\mathbb{C} ^{n + 2}$,

\item[(iii)]  the real Grassmannians
  $\widetilde{\Gr} _4(\RM^{n+4}) = \frac{\SO(n+4)}{\SO(n)\times \SO(4)}$, of
  oriented real $4$-planes in $\mathbb{R} ^{n + 4}$,
and

\item[(iv)] the five exceptional spaces
 $
  \frac{\mathrm{G}_2}{\SO(4)},\
  \frac{\mathrm{F}_4}{\Sp(3)\Sp(1)},\
  \frac{\mathrm{E}_6}{\SU(6)\Sp(1)},\
  \frac{\mathrm{E}_7}{\mathrm{Spin}(12)\Sp(1)},\
  \frac{\mathrm{E}_8}{\mathrm{E}_7\Sp(1)}
 $,
 in dimensions $4n$ with $n=2,7,10,16$ and $28$ respectively.
\end{enumerate}
According to Theorem~\ref{main}, none of them admits a (weakly)
complex structure except for the complex Grassmannians
$\Gr_2(\CM^{n+2})$. Note however that this was already known for
$\HM \P^n$, as mentioned above, and also for most real oriented
Grassmannians $\widetilde{\Gr} _4(\RM^{n+4})$. Indeed, the
non-existence of (weakly) complex structures on a large class of
real Grassmannians, including all $\widetilde{\Gr}_4(\RM^{n + 4})$
except for $\widetilde{\Gr}_4(\RM^{8})$ and
$\widetilde{\Gr}_4(\RM^{10})$, was shown in \cite{s91} by
P.~Sankaran and in \cite{z94} by Z.-Z. Tang.

Wolf spaces are (irreducible, simply-connected) {\it inner}
symmetric spaces of compact type, i.e. are symmetric spaces of the
form $G/H$, where $G, H$ are connected compact Lie groups of {\em
equal rank}. Apart from the Wolf spaces, the class of simply
connected irreducible inner symmetric spaces of compact type
includes, cf. e.g. \cite{hel78}, \cite{besse}, \cite{br90}:
\begin{enumerate}

\item[(i)]  the class of (irreducible) Hermitian
symmetric spaces of compact type~;

\item[(ii)]  the even-dimensional spheres
$\S ^{2n} = \frac{{\rm SO} (2n + 1)}{{\rm SO} (2n)}$, $n > 1$~;

\item[(iv)]  the even-dimensional oriented real Grassmannians
$\widetilde{{\rm Gr}} _{2p} (\mathbb{R} ^{n + 2p}) =
\frac{{\rm SO} (n + 2 p)}{{\rm SO}
    (n) \times {\rm SO} (2 p)}$, $n > 1$~;

\item[(v)] the quaternionic Grassmannians ${\rm Gr} _k (\mathbb{H} ^{k
    + n}) = \frac{{\rm Sp} (n + k)}{{\rm Sp} (n) \, {\rm Sp}
    (k)}$, $n, k > 1$~;

\item[(iii)]  the Cayley projective
plane $\frac{{\rm F} _4}{{\rm Spin} (9)}$;

\item[(vi)] the two exceptional inner symmetric spaces
$\frac{{\rm E} _7}{{\rm SU} (8)/ \mathbb{Z} _2}$ and $\frac{{\rm E} _8}{\Spin (16) /
  \mathbb{Z} _2}$.
\end{enumerate}
Notice that all spaces in this list are even-dimensional, cf.
Section \ref{st1}. 
In the second part of this paper, we show that the techniques
introduced in the proof of our main Theorem \ref{main} can be used to
establish a similar non-existence theorem for inner symmetric spaces
of compact type. More precisely,
we have:
\begin{ath}\label{t1} A $4n$-dimensional 
simply connected  irreducible inner symmetric space of compact
type is weakly complex if and only if it is 
a sphere or a Hermitian symmetric space.
\end{ath}
Recall that, for any $n$, the sphere $\S^n$ is stably parallelizable, hence
weakly complex, whereas
Hermitian symmetric spaces are complex manifolds in a natural way.

Our method uses a unified argument
based on the index
calculation of a twisted Dirac operator, via Theorem \ref{gregor} and
Proposition \ref{crit}
below. On the other hand, our approach is ineffective for {\it non-inner}
  symmetric spaces, due to the fact, established by R. Bott
  \cite{bott65}, that the index of
  any homogeneous differential operator vanishes on any non-inner
  symmetric space of compact type.  For
  reasons which will be explained in Section \ref{st1}, 
 cf. in particular Remarks
  \ref{r3}, \ref{r1}, \ref{r2}, it is also 
  ineffective for inner symmetric spaces of dimension $4 n + 2$. 
 In the above list, these are: (i) Hermitian symmetric
  spaces of odd complex dimension; (ii) oriented real Grassmannians 
$\widetilde{{\rm Gr}} _{2 p} (\mathbb{R} ^{2p+q})$, whith $p$ and $q$
both 
odd; (iii) the
exceptional symmetric space $\frac{{\rm E} _7}{{\rm SU} (8)/
  \mathbb{Z} _2}$ (which is of dimension $70$). 
The non-existence of weakly
complex structure for all oriented real Grassmannians, except for 
$\widetilde{{\rm Gr}} _4 (\mathbb{R} ^8)$,
$\widetilde{{\rm Gr}} _6 (\mathbb{R} ^{12})$ and
$\widetilde{{\rm Gr}} _4 (\mathbb{R} ^{10})$, in particular 
for all oriented real Grassmannians of dimensions $4 n + 2$, 
was established, by different methods, by P. Sankaran in
\cite{s91} and by Z.-Z. Tang in
\cite{z94}. 
Together with these results, our Theorem \ref{t1} then covers all
simply connected irreducible inner symmetric spaces, except for the
exceptional symmetric space $\frac{{\rm E} _7}{{\rm SU} (8)/
  \mathbb{Z} _2}$, for which, as far as we know, the existence of a
(weak) almost complex structure has remained an open question. 
Notice that the non-existence of (weak) almost complex structure
on quaternionic Grassmannians was previously established, by
a different approach,  by W. C. Hsiang
and R. H. Szczarba in \cite{hs64}. Note also 
that A. Borel and F. Hirzebruch \cite{bh58} 
have shown that the 
tangent bundle of the Cayley projective plane has no almost complex 
structure, but their proof 
does not exclude the possibility for that bundle to being weakly complex.

The irreducibility assumption in Theorem \ref{t1} can easily be
dropped. Indeed, the de Rham decomposition Theorem implies that any
simply connected inner symmetric space can be written as a product 
of irreducible
inner symmetric spaces. Using the fact that a product $M\times N$ is
weakly complex if and only if both factors are weakly complex --- the
restriction of the tangent bundle of a product to each factor being
stably isomorphic to the tangent bundle of that factor ---  and by taking 
into account the above observations, we thus
obtain the
following generalization of Theorem \ref{t1}:
\begin{ath} \label{t2} The irreducible components of a simply connected inner 
symmetric space of compact type admitting a
weak almost complex structure are isomorphic either to an even-dimensional
sphere, or to a Hermitian symmetric space or to the exceptional symmetric
space
$\frac{{\rm E} _7}{{\rm SU} (8)/ \mathbb{Z} _2}$. \end{ath} 
Our method gives however no information concerning the
existence of genuine almost complex structures on products of
even-dimensional spheres and Hermitian symmetric spaces (in contrast, 
the product $\S ^{2 p + 1} \times \S ^{2 q + 1}$ of two odd-dimensional
spheres admits {\it integrable} almost complex structures 
\cite{ce53}). Theorem \ref{t2} can be viewed as a topological
version of the well-known fact that an inner symmetric space of
compact type  which admits an {\em integrable} almost complex structure
{\em compatible
with the invariant metric}  has to be Hermitian symmetric
\cite{br90}, \cite{bmgr93}.

\noindent {\sc Acknowledgments}.  We warmly thank Claude LeBrun for
communicating the reference to Massey's paper \cite{m62} and for a
stimulating and helpful e-mail exchange. We also thank
Ulrich Bunke, Jean Lannes and Simon Salamon for useful discussions. A special
thank is due to Gregor Weingart, who informed us about 
the crucial Theorem \ref{gregor}, and
to Dieter Kotschick for comments which led to significant
improvements of the paper.

\section{Proof of Theorem \ref{main}} \label{smain}

For notation and basic properties of quaternion-K\"ahler
manifolds we refer to \cite{s82} and \cite{s89}.
Let $(M, g)$ be a $4n$-dimensional quaternion-K\"ahler manifold of
positive type, $n \geq 1$.
 Since the
holonomy group of $M$ is contained in $\Sp(n)\cdot
\Sp(1)$, the standard representations
of $\Sp(n)$ on $\CM^{2n}$ and of $\Sp(1)$ on $\CM^2$ give rise to
locally defined complex vector bundles,  denoted by $\E$ and $\H$
respectively. These are
globally defined only on the quaternionic projective space $\HM
\P^n$
(\cite[Theorem 6.3]{s82}). However, tensor products of
any {\it even}  number of copies of $\H$ and $\E$ are globally defined complex
bundles over any quaternion-K\"ahler manifold $M$.

It is well known that the complexified tangent bundle of $M$ is given
as $\T M^{\CM} = \E \otimes \H$. Recall that a
quaternion-K\"ahler manifold $M^{4n}$ of positive type
is spin if and only if either $M ^{4n} =\HM
\P^n$, or the quaternionic dimension $n$ is
even (\cite[Proposition 2.3]{s82}). If this holds, the spinor bundle $\Sigma
M$ decomposes as the direct sum of
$\R^{p,q}:=\Sym^p\H\otimes \Lambda^q_0\E$ over all positive integers $p$,
$q$ with $p+q=n$, cf. e.g. \cite[Proposition 2.1]{ksw99}. Here
$\Lambda^q_0\E$ denotes the
sub-bundle of $\Lambda^q\E$ defined as the kernel of the
contraction with the symplectic form of $\E$.  In particular, the twisted
spin bundles
$\Sigma ^\pm M \otimes \R^{p,q}$ are globally defined whenever $p+q+n$ is
even. We then denote by $\D_{\R^{p,q}}$ be the (twisted) Dirac operator defined
on sections of $\Sigma ^+ M
\otimes \R^{p,q}$ and by  $\ind(\D_{\R^{p,q}})$ the index of
$\D_{\R^{p,q}}$.

Our argument crucially relies on the following result
of C.~LeBrun and S.~Salamon
\cite[Theorem 5.1]{ls94}
(cf. also \cite{sw02}):
\begin{equation}\label{index}
\ind(\D_{\R^{p,q}}) \;=\;
\begin{cases} \quad 0 & \qquad\mbox{for} \quad p+q<n \\ \quad
  (-1)^q\,(\bb_{2q} (M) +
  \bb_{2q-2} (M)) & \qquad\mbox{for} \quad p+q=n \ ,
\end{cases}
\end{equation}
where ${\rm b} _i (M)$ denote the Betti numbers of $M$.
Consider the twist bundle $V= \Sym^{n-2}\H \otimes \T M^{\CM}$ (it is
here that the
assumption $n\ge 2$ is needed). The
Clebsch-Gordan decomposition yields
$$
V = (\Sym^{n-1}\H\otimes \E) \,\oplus\, (\Sym^{n-3}\H\otimes \E) \ .
$$
The bundle $\Sigma M \otimes V$ is globally defined for all quaternionic
dimensions $n$ and we can therefore compute the index $\ind(\D_V)$ of
the corresponding twisted Dirac operator by using
(\ref{index}). We thus obtain
\begin{equation}\label{in}
\ind(\D_{\Sym^{n-2}\H \otimes \T M^{\CM}})
\;=\;
\ind(\D_{\Sym^{n-1}\H\otimes \E}) + \ind(\D_{\Sym^{n-3}\H\otimes \E})
\;=\;
-(\bb_2 (M) + \bb_0 (M))\ .
\end{equation}

A key fact, cf. \cite[Corollary 4.3]{ls94}, is that $\bb_2(M) =
0$ for
all compact quaternion-K\"ahler
manifold $M$ of positive type other than
the complex Grassmannians $\Gr_2(\CM^{n+2})$, whereas $\bb_2 (M) = 1$
if $M = \Gr_2(\CM^{n+2})$, which, as already observed,
has a natural complex structure. We now assume
that $M$ is different from $\Gr_2(\CM^{n+2})$, so that $\bb_2 (M) =
0$.
The above index calculation then reads
\begin{equation}\label{index1}
\ind(\D_{\Sym^{n-2}\H \otimes \T M^{\CM}}) \;=\; - 1 \ .
\end{equation}
Assume, for a contradiction, that $M$ carries an almost complex
structure.
Then
the tangent bundle $\T M$ is a complex vector bundle and its
complexification splits into the sum of two complex sub-bundles
$\T M^{\CM} = \theta \oplus \theta^*$. For the components of the Chern
character we have $\ch_i(\theta^*)=(-1)^i \ch_i(\theta)$. On the
other hand, $\ch(\Sym^{n-2}\H)$ and $\hat{\mathrm{A}}(\T M)$ have non-zero
components only in degree $4k$. Indeed, $\hat{\mathrm{A}}(\T M)$ is a
polynomial in the Pontryagin classes of $M$ and $\Sym^{n-2}\H$ is a
self-dual locally defined complex bundle.
The Atiyah-Singer formula for twisted Dirac operators
(cf. \cite{bgv92}) then yields
\begin{equation}\begin{split}
\ind(\D_{\Sym^{n-2}\H \otimes \T M^{\CM}})
=&
\ch(\Sym^{n-2}\H)\ch(\T M^{\CM})\hat{\mathrm{A}}(\T M)[M]\\
=&
2\,\ch(\Sym^{n-2}\H)\ch(\theta)\hat{\mathrm{A}}(\T M)[M] \ .
\end{split}\end{equation}
Notice that $\ch(\Sym^{n-2}\H)$ is well-defined in $H ^* (M,
\mathbb{Q})$,
even if $n$ is odd.

Now, $\ch(\Sym^{n-2}\H) \ch(\theta) \hat{\mathrm{A}}(\T M)[M]$ is
the
index of the twisted Dirac operator $\D_{\Sym^{n-2}\H \otimes\theta}$
on the (globally defined) bundle $\Sigma M
\otimes\Sym^{n-2}\H \otimes\theta $ and thus has to be an
integer. This implies that
that $\ind(\D_{\Sym^{n-2}\H \otimes \T M^{\CM}})$ is {\it even}, hence
contradicts (\ref{index1}).

If the manifold is assumed to be weakly complex then
there exists a trivial real vector bundle $\e$ such that $\T M \oplus \e$
is a complex vector bundle. By replacing $V= \Sym^{n-2}\H \otimes \T
M^{\CM}$ with $V=\Sym^{n-2}\H \otimes
(\T M \oplus \e)^{\CM}$ in the above argument, this
remains unchanged, as the extra term
$$
\ind(\D_{\Sym^{n-2}\H \otimes \e^{\CM}})
=
\rk(\e)\,\ind(\D_{\Sym^{n-2}\H})
$$
in \eqref{in} is zero, again because of (\ref{index}). This completes the
proof of Theorem~\ref{main}.

\section{Proof of Theorem \ref{t1}}\label{st1}

We first establish a
general formula, of separate interest,  for the index
of a family of homogeneous twisted Dirac operators defined on inner
symmetric spaces.

Let $M=G/K$ be an irreducible inner symmetric space of compact type,
where $G$ denotes a (connected) compact simple Lie group and $K$ a
connected closed subgroup of $G$.
Notice that the condition implies that $M$ is even-dimensional.
We fix a
common maximal torus $T\subset K\subset G$ and we endow the dual Lie
algebra $\tt^*$ with a suitable positive definite
scalar product $\langle\cdot,\cdot\rangle$,
proportional to the one induced by the opposite of the Killing form of
$G$.
We denote by $\rho^\kk$ and $\rho^\gg$
the half-sum of the positive roots of $K$ and $G$ respectively.

The isotropy
representation $K\to\SO(\mm)$ induces a group homomorphism $\tilde K\to
\Spin(\mm)$ and thus a
representation of $\tilde K$ on the spin modules $\Sigma^\pm_\mm$,
where $\tilde{K}$ stands for $K$ itself or a two-fold covering of $K$.
Let $V _{\mu}$ be a complex representation of $\tilde{K}$ with highest
weight $\mu \in \mathfrak{t} ^*$. We assume
that the induced representation  of $\tilde K$ on $V_\mu\otimes
\Sigma^\pm_\mm$ descends to a representation of $K$. We then denote by
$\Sigma^\pm_\mu M:=G\times_K (V_\mu\otimes
\Sigma^\pm_\mm)$ the corresponding
{\it twisted spin bundles}   and by $\D _{V_\mu}$ the
{\it twisted Dirac operator}
acting on sections of $\Sigma^\pm_\mu M:=G\times_K (V_\mu\otimes
\Sigma^\pm_\mm)$.

\begin{ath}\label{gregor}
Let $w \in W_\gg$ be a Weyl group element for which $w(\mu +
\rho^\kk) - \rho^\gg$ is $\gg$-dominant. Then the index of the
twisted Dirac operator $\D _{V_\mu}:C^\infty(\Sigma^+_\mu M)\to
C^\infty(\Sigma^-_\mu M)$ is given by the formula \be\label{im}
  \ind(\D _{V_\mu})
  \;\; = \;\;
  \prod_{\alpha\in\RR^+}
  \frac{\langle \mu + \rho^\kk, \alpha \rangle}{\langle
    \rho^\gg,\alpha \rangle} =: \i(\mu),
\ee where the product goes over all positive roots of $G$. If such a
$w$ does not exist the index is zero.
\end{ath}
\begin{proof}
The generalized Bott-Borel-Weil theorem (cf. \cite{huang06},
Theorem~4.5.1) states that the kernel of $\D _{V_\mu}$ is an
irreducible $G$-representation of highest weight $w(\mu + \rho^\kk)
- \rho^\gg$, where $w$ is as in the assumption of the theorem. Note
that $w$ is unique as soon as it exists. If such a $w$ does not
exist, the kernel is zero. It follows that the index of $\D _{V_\mu}$
is given as $(-1)^{l(w)}\dim V_{w(\mu + \rho^\kk) - \rho^\gg}$, or
zero if such a $w$ does not exist (cf. \cite{huang06},
Corollary~4.5.2). Here $l(w)$ is the length of the Weyl group
element $w$, defined as the number of positive roots $\alpha$ such
that $w(\alpha)$ is a negative root. It is also the smallest integer
$k$ for which $w$ can be written as the product of $k$ reflections
in simple roots. Thus the length of $w$ is the same as the length of
$w^{-1}=w^t$. The Weyl dimension formula implies
\begin{equation}\label{dim}
\dim (V_{w(\mu + \rho^\kk) - \rho^\gg}) = \prod_{\alpha
\in\RR^+}\frac{\langle w(\mu + \rho^\kk), \alpha \rangle}{\langle
    \rho^\gg,\alpha \rangle}
    =
    \prod_{\alpha \in\RR^+} \langle \rho^\gg,\alpha \rangle^{-1}
    \prod_{\alpha \in\RR^+}\langle \mu +
\rho^\kk, w^t(\alpha) \rangle
\end{equation}
If we replace in the set $\{w^t(\alpha) |\;  \alpha \in \RR^+\}$ the
$l(w)$ roots which are mapped by $w^t$ to negative roots by their
negative, we obtain again the set $\RR^+$ of positive roots. Hence
$$
\prod_{\alpha \in\RR^+}\langle \mu + \rho^\kk, w^t(\alpha) \rangle =
(-1)^{l(w)}\prod_{\alpha \in\RR^+}\langle \mu + \rho^\kk, \alpha
\rangle.
$$
Substituting this into formula (\ref{dim}) and using
$\ind(\D _{V_\mu})=(-1)^{l(w)}\dim V_{w(\mu + \rho^\kk) - \rho^\gg}$
completes the proof of the theorem.
\end{proof}

The following criterion, extracted from the proof of Theorem \ref{main}
in Section \ref{smain}, provides a general obstruction for the
tangent bundle of a compact manifold to being weakly complex.

\begin{epr} \label{crit} Let $(M^{4n},g)$ be a compact Riemannian
  manifold carrying a locally defined complex vector
  bundle $E$ such that the following conditions hold:
\begin{enumerate}
\item[(a)] $E$ is self-dual, i.e. $E$ is isomorphic to its dual bundle $E^*$.
\item[(b)] $E\otimes \Sigma M$ is globally defined, where $\Sigma M$
  denotes the (locally defined) spin bundle of $M$.
\item[(c)] The index of the twisted Dirac operator
  $\D_{E\otimes\T M^\CM}$ is odd.
\end{enumerate}
Then the tangent bundle of $M$ is not weakly complex and in particular
$M$ cannot
carry an almost complex structure.
\end{epr}

\begin{proof}
The Atiyah-Singer index formula for twisted Dirac operators
(cf. \cite{bgv92}) reads
\begin{equation}\label{as}
\ind(\D_{V})=
\ch(V)\hat{\mathrm{A}}(\T M)[M]\end{equation}
for every complex bundle $V$ such that $V\otimes \Sigma M$ is globally defined.

Assume that $\T M$ is weakly complex, i.e. there exists a trivial real
vector bundle $\e$ (of even rank) such that $\T M \oplus \e$ is a
complex vector bundle.
Then its  complexification splits into the sum of two complex bundles
$(\T M\oplus \e)^{\CM} = \theta \oplus \theta^*$.

Since $\hat{\mathrm{A}}(\T M)$ is a
polynomial in the Pontryagin classes of $M$, it has non-zero components
only in degree $4k$.
Condition (a) shows that the Chern character $\ch(E)$ has the same
property. Moreover, the components of the Chern
characters of $\theta$ and $\theta^*$ satisfy $\ch_i(\theta^*)=(-1)^i
\ch_i(\theta)$. Applying \eqref{as} to $V=E\otimes \theta$ and
$V=E\otimes \theta^*$ yields
\be \ind(\D_{E\otimes\theta})=\ch(E)\ch(\theta)\hat{\mathrm{A}}(\T
M)[M]=\ch(E)\ch(\theta^*)\hat{\mathrm{A}}(\T
M)[M]=\ind(\D_{E\otimes\theta^*}).
\ee
Using this equation, the condition that $\rank(\e)$ is even and assumption (b), we infer
\bea \ind(\D_{E\otimes\T M^\CM})&\equiv&\ind(\D_{E\otimes\T
  M^\CM})+\rank(\e)\ind(\D_{E})\ \mod 2\\
&\equiv&\ind(\D_{E\otimes(\T M^\CM\oplus \e^{\CM})})\ \mod 2\\
&\equiv&\ind(\D_{E\otimes\theta}) + \ind(\D_{E\otimes\theta^*})\ \mod 2\\
&\equiv&2\ind(\D_{E\otimes\theta})\equiv 0\ \mod 2,
\eea
contradicting (c). This proves the proposition.
\end{proof}

\begin{ere}\label{r3}
Notice that in dimension $4n+2$, there exist no local bundle $E$ satisfying conditions (a)--(c). Indeed, if $E$ is self-dual, the Chern character of $E\otimes \T M^{\CM}$ has non-zero components only in degree $4k$. Formula \eqref{as} applied to $V=E\otimes \T M^{\CM}$ then shows that the index of the twisted Dirac operator $\D_{E\otimes\T M^\CM}$ vanishes. 
\end{ere}

We now check that the criterion given by Proposition \ref{crit}
applies to all $4n$-dimensional (simply connected) irreducible
inner symmetric spaces of compact type, using as main tool
the formula (\ref{im}) in  Theorem \ref{gregor}. We focus on those cases which were
not fully covered by previous works, namely the oriented real
Grassmannians $\widetilde \Gr_{2p}(\RM^{2p+q})$ 
with either $p$ or $q$ even (of dimension $2pq$), the exceptional
inner symmetric space $\E_8/(\Spin(16)/\ZM_2)$ (of dimension 128), and the Cayley projective plane
$\F_4/\Spin(9)$ (of dimension 16). 
The quaternionic Grassmannians can be handled with quite similar methods.

\subsection{The oriented real
Grassmannians $\widetilde \Gr_{2p}(\RM^{2p+q})$} 

\subsubsection{Case I: $q=2q'$ is even.}
Since for $p=1$ or $q'=1$ the Grassmannian of oriented 2-planes is a Hermitian symmetric space, we
assume $p,q'\ge 2$. The symmetric space $M=G/K:=\SO(2p+2q')/\SO(2p)\times \SO(2q')$ is spin
(cf. \cite{hs90}). Let $H$ and $H'$ denote the
tautological bundles over $M$, associated to the standard
representations of $\SO(2p)\times \SO(2q')$ on $\RM^{2p}$ and
$\RM^{2q'}$ respectively. It is well-known that $\T M$ is isomorphic to
$H\otimes H'$ (cf. \cite{besse}, p. 312).

The root system of $G$ consists of the vectors $\pm e_i\pm e_j$,
$1\le i< j\le p+q'$. We choose as fundamental Weyl chamber the one
containing the vector $\sum_{i=1}^{p+q'}(p+q'-i)e_i$. The positive roots
are then $e_i\pm e_j$, $1\le i< j\le p+q'$, and their half-sum is
$$\rho ^\gg=\sum_{i=1}^{p+q'}(p+q'-i)e_i.$$
The root system of $K$ is the direct sum of the root systems of
$\SO(2p)$ and $\SO(2q')$ and thus
$$\rho ^\kk=\sum_{i=1}^{p}(p-i)e_i+\sum_{j=1}^{q'}(q'-j)e_{p+j}.$$
Let $E$ be the complex vector bundle over $M$ associated to the
complex representation of $K$ with highest weight
\beq\label{a1}\mu=(q'-1)\sum_{i=1}^{p}e_i+(p-2)e_{p+1}.\eeq
In other words, $E$ is the Cartan component in the tensor product of 
the $(q'-1)$-th symmetric power of $\Lambda ^p H^\CM$ and the $(p-2)$-th symmetric power of
$(H')^\CM$. It clearly satisfies (a) and (b) in Proposition \ref{crit}, and
we claim that it also satisfies (c).

To see this, we need to compute the decomposition in irreducible
summands of $E\otimes \T M^\CM$. This is given by the following
standard facts:

\begin{elem}\label{l1}
The tensor product of complex irreducible $\Spin(2r)$ representations with
highest weights $e_1$ and $k(e_1+\ldots +e_r)$ is the direct sum of
representations with highest weights $(k+1)e_1+k(e_2+\ldots +e_r)$ and
$k(e_1+\ldots +e_{r-1})+(k-1)e_r$.

The tensor product of complex irreducible $\Spin(2r)$ representations with
highest weights $e_1$ and $ke_1$ is the direct sum of representations with
highest weights $(k+1)e_1$, $ke_1+e_2$ and $(k-1)e_1$.\end{elem}

\begin{proof}
The statements of the lemma follow from a routine decomposition of
tensor products. However, in this special case, where one of the
factors is the standard representation of   $\Spin(2r)$,  the
decomposition can directly be read off the table given in
\cite{gu10} p. 511.
\end{proof}

The complexified bundles $H^\CM$ and $(H')^\CM$ are
associated to the irreducible $\SO(2p)\times \SO(2q')$
representations with highest weights $e_1$ and $e_{p+1}$.
(This is where the hypothesis $p,q'\ge 2$ is used:
The complexification of the standard representation
of $\SO(2)$ on $\RM^2$ is reducible!)
Lemma \ref{l1} shows that $E\otimes \T M^\CM$ is associated to
the direct sum of representations with highest weights
$$
\mu_1=q'e_1+(q'-1)\sum_{i=2}^{p}e_i+(p-1)e_{p+1},\qquad\ \
\mu_2=(q'-1)\sum_{i=1}^{p-1}e_i+(q'-2)e_{p}+(p-1)e_{p+1},\ \ \ \ $$
$$\mu_3=q'e_1+(q'-1)\sum_{i=2}^{p}e_i+(p-2)e_{p+1}+e_{p+2},\quad
\mu_4=(q'-1)\sum_{i=1}^{p-1}e_i+(q'-2)e_{p}+(p-2)e_{p+1}+e_{p+2},$$
$$\mu_5=q'e_1+(q' - 1)\sum_{i=2}^{p}e_i+(p-3)e_{p+1},\qquad\ \
\mu_6=(q'-1)\sum_{i=1}^{p-1}e_i+(q'-2)e_{p}+(p-3)e_{p+1}.\ \ \ \
$$
It is clear that the coordinates of $\mu_1+\rho ^\kk$ are a
permutation of the coordinates of $\rho ^\gg$, so
$i(\mu_1)=\pm1$. Moreover, $i(\mu_i)=0$ for $2\le i\le 6$ since
$\mu_i+\rho ^\kk$ has two equal coordinates in each case. By Theorem
\ref{gregor}, condition (c) in Proposition \ref{crit} is satisfied, hence
  the tangent bundle of $M$ is not weakly complex.

\subsubsection{Case II: $p=2p'$ is even.}
We can assume that $q$ is odd since the case 
where $q$ is even is included in the previous one. 
The manifold $M=\SO(2p+q)/\SO(2p)\times
\SO(q)$ is not spin (cf. \cite{hs90}). Nevertheless, the tensor product of
the (locally defined) spin bundle of $M$ and any 
odd tensor power of the locally
defined spin bundle of $H$ is globally defined. For $q\ge 3$ we take $E$ as
the locally defined complex vector bundle over $M$ associated to the
complex representation of $\Spin(2p)\times\Spin(q)$ with highest weight
\beq\label{a2}\mu=(\tfrac q2-1)\sum_{i=1}^{p}e_i+(p-2)e_{p+1}\eeq
(incidentally this is exactly the same formula as \eqref{a1}).
Like before, Theorem \ref{gregor} shows that the index of the Dirac operator twisted with $E\otimes \T M^{\CM}$ is $\pm 1$. Moreover, $E$ is self-dual, being associated to the Cartan component of the 
tensor product of self-dual representations: The complexification of the standard representation of $\SO(q)$ and the spin representation of $\Spin(4p')$. Proposition \ref{crit} thus shows that $M$ is not weakly complex.

The argument does not apply for $q=1$ since $\mu$ is no
longer a highest weight in that case. This is of course compatible
with the fact that the sphere $S^{2p}= \SO(2p+1)/\SO(2p)$ is weakly complex,
having stably trivial tangent bundle.

\begin{ere}\label{r1}
The above construction can also be carried out verbatim on the remaining even-dimensional oriented real
Grassmannians $\widetilde \Gr_{2p}(\RM^{2p+q})$, when  $p$ and $q$ are
both odd, by choosing for $E$ 
the locally defined complex vector bundle associated to the
complex representation of $\Spin(2p)\times\Spin(q)$ with highest weight
given by \eqref{a2}.
In other words, $E$ is the Cartan component in the tensor product of 
the $(q-2)$-th symmetric power of the spin representation $\Sigma^+_{2p}$ 
and the $(p-2)$-th symmetric power of
$(H')^\CM$.
The same argument shows that the index of the corresponding 
twisted Dirac operator 
is $\pm1$. However, the bundle $E$ is no longer self-dual 
if $p$ and $q$ are both odd, since
the spin representation of $\Sigma^+_{2p}$ -- and 
thus its $(q-2)$-th symmetric power -- is not self-dual in this case.
\end{ere}

\subsection{The exceptional symmetric space $M=\E_8/(\Spin(16)/\ZM_2)$}
The group $\ZM_2$ acting on $\Spin(16)$ is generated by
the volume element $v:=e_1\,\ldots\,e_{16}\in \Spin(16)$
(cf. \cite{a}). The positive
half-spin representation factors through $v$ and induces a
representation of $\Spin(16)/\ZM_2$ on $\Sigma^+_{16}$ whose
associated bundle is just $\T M^\CM$. Since $v$ maps to $-\id\in
\SO(16)$, the representation of $\Spin(16)$ on $\RM^{16}$ defined by
the spin covering $\xi:\Spin(16)\to \SO(16)$ induces a locally defined
real vector bundle $H$ on $M$. Of course, all even tensor products of $H$ are
globally defined vector bundles on $M$. Moreover, the manifold $M$ is
spin (cf \cite{hs90}). Conditions (a) and (b) in Proposition \ref{crit} are
thus satisfied for $E=\Sym^{2k}_0H^\CM$, i.e. the Cartan summand of
the $2k$-th symmetric power of $H^\CM$, associated to the
representation of $\Spin(16)$ with highest weight
$(2k,0,0,0,0,0,0,0)\in \tt^*\simeq \RM^8$.

We will use Theorem \ref{gregor} in order to compute the index of the
Dirac operator on $M$ twisted with $\T M^\CM\otimes
\Sym^{2k}_0H^\CM$. Since $\T M^\CM$ is associated to the positive half-spin
representation, whose highest weight is $\tfrac12(1,1,1,1,1,1,1,1)$,
we need to decompose $\Sym^{2k}_0H^\CM\otimes \Sigma_{16}^+$ into irreducible
components.

\begin{elem}\label{dec} $\Sym^{2k}_0H^\CM\otimes \Sigma_{16}^+$ is the
direct sum of the $\Spin(16)$-representations with highest weights
$\tfrac12(4k+1,1,1,1,1,1,1,1)$ and $\tfrac12(4k-1,1,1,1,1,1,1,-1)$.
\end{elem}
\begin{proof}
Again the decomposition follows from a standard calculation, where
in this case the result can also be found in \cite{ov90}, p. 303.
\end{proof}

The root system $\RR(\E_8)$ is the disjoint union of the root
system of $\Spin(16)$ and the weights of the half-spin representation
$\Sigma^+_{16}$. It thus consists of the vectors $\pm e_i\pm e_j$,
$1\le i< j\le 8$ and
$$\frac 12 \sum_{i=1}^8 \e_i e_i,\qquad \e_i=\pm 1,\ \e_1\cdots
\e_8=1.$$
With respect to the fundamental Weyl chamber containing the vector
$(23,6,5,4,3,2,1,0)$, the set of positive roots of $\Spin(16)$ is
$\RR^+(\Spin(16))=\{e_i\pm e_j\ |\ 1\le i< j\le 8\}$, and the set of
positive roots of $\E_8$ is
$$\RR^+(\E_8)=\RR^+(\Spin(16))\cup\bigg\{\frac 12 \sum_{i=1}^8 \e_i e_i\ |\
\e_1=1,\ \e_i=\pm 1,\ \e_1\cdots \e_8=1\bigg\}.$$
The half-sums of the positive roots of $K=\Spin(16)/\ZM_2$ and
$G=\E_8$ are thus given by
$$\rho^\kk=(7,6,5,4,3,2,1,0)\quad \hbox{and}\quad
\rho^\gg=(23,6,5,4,3,2,1,0).$$
It is clear that $\mu_1+\rho
  ^\kk$ is orthogonal to the root
  $\a=\tfrac12(1,-1,-1,-1,-1,-1,-1,1)$ for
  $\mu_1:=\tfrac12(33,1,1,1,1,1,1,1)$, so the integer $\i(\mu)$ defined
  by \eqref{im}
 vanishes for $\mu=\mu_1$. Moreover, an elementary
 computation shows
 that $\i(\mu_2)=-1$ for
  $\mu_2:=\tfrac12(31,1,1,1,1,1,1,-1)$. By Lemma \ref{dec}, the
  tensor product
  $\T M^\CM\otimes \Sym^{16}_0H^\CM$ is associated to the direct sum of
  $\Spin(16)$-representations with highest weights $\mu_1$ and
  $\mu_2$. Theorem \ref{gregor} thus shows that condition (c) in
  Proposition \ref{crit} is satisfied for $E=\Sym^{16}_0H^\CM$, hence
  the tangent bundle of $M$ is not weakly complex.

\subsection{The Cayley projective plane $M=\F_4/\Spin(9)$}
The complexified tangent bundle $\T M^\CM$ is associated to the spin 
representation on $\Sigma_9\simeq\CM^{16}$ (cf. \cite{besse}, p. 302). 
Let $H$ denote the real bundle associated to the representation 
of $\Spin(9)$ on $\RM^9$ defined by the spin covering $\Spin(9)\to\SO(9)$,
with highest weight
$(1,0,0,0)\in \tt^*\simeq \RM^4$.
It is well-known that the Cayley projective plane is
spin (cf \cite{hs90}). Conditions (a) and (b) in Proposition \ref{crit} are
thus satisfied for $E=H^\CM$. 

We will use Theorem \ref{gregor} again in order to compute the index of the
Dirac operator on $M$ twisted with $\T M^\CM\otimes H^\CM$. Recall that
$\Sigma_9\otimes \CM^9\simeq \Sigma_9\oplus\Sigma_9^{\frac32}$, the 
two summands having highest weights $\tfrac12(1,1,1,1)$ and $\tfrac12(3,1,1,1)$
respectively.

The root system $\RR(\F_4)$ is the disjoint union of the root
system of $\Spin(9)$ and the weights of the spin representation
$\Sigma_{9}$. It thus consists of the vectors $\pm e_i\pm e_j$, 
$1\le i< j\le 4$, $\pm e_i$, $1\le i\le 4$ and
$$\frac 12 \sum_{i=1}^4 \e_i e_i,\qquad \e_i=\pm 1.$$
With respect to the fundamental Weyl chamber containing the vector
$(11,5,3,1)$, the set of positive roots of $\Spin(9)$ is
$\RR^+(\Spin(9))=\{e_i\pm e_j\ |\ 1\le i< j\le 4\}\cup \{e_i\ |\ 1\le i\le 4\}$, and the set of
positive roots of $\F_4$ is
$$\RR^+(\F_4)=\RR^+(\Spin(9))\cup\bigg\{\frac 12 \sum_{i=1}^4 \e_i e_i\ |\
\e_1=1,\ \e_i=\pm 1\bigg\}.$$
The half-sums of the positive roots of $K=\Spin(9)$ and
$G=\F_4$ are thus given by
$$\rho^\kk=\tfrac12(7,5,3,1)\quad \hbox{and}\quad
\rho^\gg=\tfrac12(11,5,3,1).$$
It is clear that $\tfrac12(1,1,1,1)+\rho
  ^\kk$ is orthogonal to the root
  $\a=\tfrac12(1,-1,-1,1)$ so the integer $\i(\mu)$ defined
  by \eqref{im}
 vanishes for $\mu=\tfrac12(1,1,1,1)$. An easy elementary
 computation shows that $\i(\mu)=-1$ for
  $\mu:=\tfrac12(3,1,1,1)$. Theorem \ref{gregor} thus shows that condition (c) in
  Proposition \ref{crit} is satisfied for $E=H^\CM$, hence
  the tangent bundle of $M$ is not weakly complex.

\medskip

\begin{ere}\label{r2}
A similar argument shows that the remaining exceptional symmetric space $M=\E_7/(\SU(8)/\ZM_2)$ (of dimension 70) also carries a complex vector bundle $E$ satisfying conditions (b) and (c)
in Proposition \ref{crit}. More precisely, $E$ is the Cartan component of the
10-th symmetric power of the locally defined bundle on $M$ associated
to the standard representation of $\SU(8)$. Of course, $E$ is not self-dual, cf. Remark \ref{r3}. 
\end{ere}

\end{document}